\documentclass[12pt,twoside,a4paper]{amsart}
\usepackage{amssymb,verbatim,enumerate,ifthen,hyperref}
\usepackage[utf8]{inputenc}
\usepackage[T1]{fontenc}
\usepackage[mathscr]{eucal}

\newtheorem*{theoremG}{Theorem G}
\newtheorem*{theoremS}{Theorem S} 
\newtheorem{lemma}{Lemma}
\newcommand{\Lem}[2]{\begin{lemma}\label{L#1}#2\end{lemma}}
\newcommand{\lem}[1]{Lemma \ref{L#1}}

\newtheorem*{corollaryG}{Corollary G}
\newtheorem*{corollaryS}{Corollary S}

\theoremstyle{remark}

\def\eq#1{{\rm(\ref{E#1})}}
\def\Eq#1#2{\ifthenelse{\equal{#1}{*}}
  {\begin{equation*}\begin{aligned}#2\end{aligned}\end{equation*}}
  {\begin{equation}\begin{aligned}\label{E#1}#2\end{aligned}\end{equation}}}

\newcommand{\R}{\mathbb{R}}
\newcommand{\N}{\mathbb{N}}

\newcommand{\A}{\mathscr{A}}
\newcommand{\G}{\mathscr{G}}
\renewcommand{\H}{\mathscr{H}}

\newcommand{\AG}{\mathscr{A}\otimes\mathscr{G}}

\begin{document}
\begin{flushright}
Appl. Math. Comput \textbf{216} (2010), 3219–3227. \\
\href{http://dx.doi.org/10.1016/j.amc.2010.04.046}{doi: 10.1016/j.amc.2010.04.046} \\[1cm]
\end{flushright}

\title[Invariance equation for two-variable Stolarsky means]
{Computer aided solution of the invariance equation for two-variable Stolarsky means}
\author[Sz. Baják]{Szabolcs Baják}
\email{bajaksz@gmail.com}
\address{Faculty of Informatics, University of Debrecen, 
4010 Debrecen, Pf 12, Hungary}
\author[Zs. Páles]{Zsolt Páles}
\email{pales@science.unideb.hu}
\address{Institute of Mathematics, University of Debrecen, 
4010 Debrecen, Pf 12, Hungary}
\date\today
\thanks{\textit{Corresponding author:} Zsolt Páles, 
Institute of Mathematics, University of Debrecen, 
4010 Debrecen, Pf 12, Hungary}
\thanks{This research was supported by the Hungarian Research Fund (OTKA)
  Grant Nos. NK-68040, NK-81402.}

\begin{abstract}
We solve the so-called invariance equation in the class of two-variable 
Stolarsky means $\{S_{p,q}:p,q\in\R\}$, i.e., we find necessary and sufficient 
conditions on the 6 parameters $a,b,c,d,p,q$ such that the identity
\[
   S_{p,q}\big(S_{a,b}(x,y),S_{c,d}(x,y)\big)=S_{p,q}(x,y) \qquad (x,y \in \R_+)
\]
be valid. We recall that, for $pq(p-q)\neq 0$ and $x\neq y$, 
the Stolarsky mean $S_{p,q}$ is defined by
\[
  S_{p,q}(x,y):=\left(\dfrac{q(x^p-y^p)}{p(x^q-y^q)}\right)^{\frac1{p-q}}.
\]
In the proof first we approximate the Stolarsky mean and we use the computer algebra 
system Maple V Release 9 to compute the Taylor expansion of the approximation
up to 12th order, which enables us to describe all the cases of the equality. 
\end{abstract}

\maketitle

\section{Introduction}

Let $\R_+$ denote the set of positive real numbers
throughout this paper. A two-variable continuous
function $M:\R_+^2\to \R_+$ is called a \emph{mean} on $\R_+$ if  
\Eq{M}{
  \min(x,y)\leq M(x,y)\leq\max(x,y)\qquad(x,y\in \R_+)
}
holds. If both inequalities in \eq{M} are strict whenever $x\neq y$,
then $M$ is called a \emph{strict mean} on $\R_+$.

Given three strict means $M,N,K:\R_+^2\to \R_+$, we say that the 
triple $(M,N,K)$ satisfies the \emph{invariance equation} if
\Eq{MNK}{
  K\big(M(x,y),N(x,y)\big)=K(x,y)\qquad(x,y\in \R_+)
}
holds.
If \eq{MNK} holds then we say that \emph{$K$ is invariant with
respect to the mean-type mapping $(M,N)$}. 
It is well known that $K$ is uniquely determined by $M$ and $N$, 
and it is called the Gauss composition $K=M\otimes N$ of $M$ and $N$.
For this terminology and result, see for instance the papers
\cite{Mat99a}, \cite{DarPal02c}.

The simplest example when the invariance equation holds is the
well-known identity
\Eq{*}{
  \sqrt{xy}=\sqrt{\dfrac{x+y}2\cdot\dfrac{2xy}{x+y}}\qquad(x,y\in\R_+),
}
that is,
\Eq{*}{
  \G(x,y)=\G\big(\A(x,y),\H(x,y)\big)\qquad(x,y\in\R_+),
}
where $\A,\G$, and $\H$ stand for the two-variable arithmetic,
geometric, and harmonic means, respectively. Another less trivial
invariance equation is the identity 
\Eq{*}{
  \AG(x,y)=\AG\big(\A(x,y),\G(x,y)\big)\qquad(x,y\in\R_+),
}
where $\AG$ denotes Gauss' \textit{arithmetic-geometric mean} defined by
\Eq{*}{
  \AG(x,y)=\bigg( \frac2{\pi}\int\limits_0^{\frac{\pi}2}
   \frac{dt}{\sqrt{x^2\cos^2t+y^2\sin^2 t}} \bigg)^{-1}\qquad
   (x,y\in\R_+). 
}
The reader is recommended to consult the book \cite{BorBor87} for more details 
and history of this deep theory.

The invariance equation in more general classes of means has recently been 
studied extensively by many authors in various papers. The invariance of the arithmetic 
mean $\A$ (i.e., when in \eq{MNK} $K$ is the arithmetic mean)
with respect to two quasi-arithmetic means was first investigated by Matkowski
\cite{Mat99a} under twice continuous differentiability assumptions concerning
the generating functions of the quasi-arithmetic means. These regularity
assumptions were weakened step-by-step by Dar\'oczy, Maksa, and Páles in the
papers \cite{DarMakPal00}, \cite{DarPal01a}, and finally this problem was completely 
solved assuming only continuity of the unknown functions involved \cite{DarPal02c}.
The invariance equation involving three weighted quasi-arithmetic means was
studied by Burai \cite{Bur07} and Jarczyk--Mat\-kowski
\cite{JarMat06}, Jarczyk \cite{Jar07}. The final answer (where no additional
regularity assumptions are required) has been obtained in \cite{Jar07}.
In a recent paper, we have studied the invariance of the arithmetic mean with 
respect to two so-called generalized quasi-arithmetic mean under four times
continuous differentiability assumptions \cite{BajPal09a}.
The invariance of the arithmetic mean with respect to Lagrangian means was the
subject of investigation of the paper \cite{Mat05} by Matkowski.
The invariance of the arithmetic, geometric, and harmonic means with respect
to the so-called Beckenbach--Gini means was studied by Matkowski in \cite{Mat02a}.
Pairs of Stolarsky means for which the geometric mean is invariant were
determined by B{\l}asi{\'n}ska-Lesk--G{\l}a\-zowska--Matkowski
\cite{BlaGlaMat03}. The invariance of the arithmetic mean with respect to
further means was studied by G{\l}azowska--Jarczyk--Matkowski
\cite{GlaJarMat02}, Burai \cite{Bur06} and Domsta--Matkowski \cite{DomMat06}.

An important class of two variable homogeneous means are the so-called
Stolarsky means (cf.\ Stolarsky \cite{Sto75}). Given two parameters $p,q\in\R$,
the two-variable mean $S_{p,q}:\R_+^2\to \R_+$ is defined by the 
following formula
\Eq{*}{
  S_{p,q}(x,y):=
  \begin{cases}
  	\bigskip
  	\left(\dfrac{q(x^p-y^p)}{p(x^q-y^q)}\right)^{\frac1{p-q}}, 
             & \mbox{ if } (p-q)pq\neq 0,\, x\neq y\\
  	\exp \left(-\dfrac1p + \dfrac{x^p \log x - y^p \log y}{x^p - y^p}\right), 
             & \mbox{ if } p=q\neq 0,\, x\neq y,\\
  	\left(\dfrac{x^p-y^p}{p(\log x - \log y)}\right)^{\frac{1}{p}}, 
             & \mbox{ if } p\neq 0,\, q=0,\, x\neq y,\\
        \left(\dfrac{q(\log x - \log y)}{x^q-y^q}\right)^{-\frac{1}{q}}, 
             & \mbox{ if } p=0,\, q\neq 0,\, x\neq y,\\
  	\sqrt{xy}, & \mbox{ if } p=q=0,\\
  	x, & \mbox{ if } x=y,
  \end{cases}
}
for every $x,y \in \R_+$.

It can easily be seen that the power (or H\"older) mean of exponent $p$
can be obtained as $S_{2p,p}$. In particular, $S_{2,1},\, S_{0,0}$ and 
$S_{-2,-1}$ are the arithmetic, geometric, and harmonic means, respectively.
Moreover, it follows by a simple computation that $S_{p,q}$ is equal to 
the geometric mean whenever $p+q=0$.
One can also check that the Stolarsky means are symmetric with respect to their
variables and their parameters as well, i.e., $S_{p,q}(x,y)=S_{p,q}(y,x)=S_{q,p}(x,y)$
for all $x,y\in \R_+$ and $p,q\in\R$.

The aim of this paper is to solve the invariance equation in the class
of Stolarsky means, i.e., to solve \eq{MNK} when each of the means $M,N$
and $K$ is a Stolarsky mean. More precisely, we want to describe the set of all 
6-tuples $(a,b,c,d,p,q)$ such that the identity
\Eq{IE}{
  S_{p,q}\big(S_{a,b}(x,y),S_{c,d}(x,y)\big)=S_{p,q}(x,y)
    \qquad(x,y\in\R_+)
}
holds. The analogous problem concerning the so-called Gini means, 
i.e., the description of the set of all 6-tuples $(a,b,c,d,p,q)$ such that
\Eq{IEG}{
  G_{p,q}\big(G_{a,b}(x,y),G_{c,d}(x,y)\big)=G_{p,q}(x,y)
    \qquad(x,y\in\R_+)
}
be valid, has recently been solved
in the paper \cite{BajPal09b} by the authors. Recall that given two parameters $p,q\in\R$,
the two-variable mean $G_{p,q}:\R_+^2\to \R_+$ is defined by the following formula 
(cf.\ Gini \cite{Gin38}):
\Eq{*}{
  G_{p,q}(x,y)=
  \begin{cases}
  	\bigskip
  	\left(\dfrac{x^p+y^p}{x^q+y^q}\right)^{\frac1{p-q}} 
           & \mbox{ for } p\neq q,\\
  	\exp \left({\dfrac{x^p \log x + y^p \log y}{x^p+y^p}}\right) 
           & \mbox{ for } p=q,
  \end{cases}
}
for $x,y \in\R_+$. The class of Gini means is also a generalization of the class of power means,
since taking $q=0$, we immediately get the power (or H\"{o}lder) mean of exponent $p$.
The main result of the paper \cite{BajPal09b} is contained in the following theorem:

\begin{theoremG}
Let $a,b,c,d,p,q\in\R$. Then the invariance equation \eq{IEG}
is satisfied if and only if one of the following possibilities hold:
\begin{enumerate}[(i)]
\item $a+b=c+d=p+q=0$, i.e., all the three means are equal to the geometric mean,
\item $\{a,b\}=\{c,d\}=\{p,q\}$, i.e., all the three means are equal to each other,
\item $\{a,b\}=\{-c,-d\}$ and $p+q=0$, i.e., $G_{p,q}$ is the geometric mean and 
      $G_{a,b}=G_{-c,-d}$,
\item there exist $u,v\in\R$ such that $\{a,b\}=\{u+v,v\}$, $\{c,d\}=\{u-v,-v\}$, and 
      $\{p,q\}=\{u,0\}$ (in this case, $G_{p,q}$ is a power mean),
\item there exists $w\in\R$ such that $\{a,b\}=\{3w,w\}$, $c+d=0$, and $\{p,q\}=\{2w,0\}$
      (in this case, $G_{p,q}$ is a power mean and $G_{c,d}$ is the geometric mean),
\item there exists $w\in\R$ such that $a+b=0$, $\{c,d\}=\{3w,w\}$, and $\{p,q\}=\{2w,0\}$
      (in this case, $G_{p,q}$ is a power mean and $G_{a,b}$ is the geometric mean).
\end{enumerate}
\end{theoremG}

As an obvious consequence, we obtain the following solution for the so-called 
Matkowski--Sut\^o equation, i.e.,
when $G_{p,q}$ is equal to the arithmetic mean in the invariance equation \eq{IEG},
which happens exactly when $\{p,q\}=\{1,0\}$.

\begin{corollaryG}
Let $a,b,c,d\in\R$. Then the Matkowski--Sut\^o equation 
\Eq{*}{
  G_{a,b}(x,y)+G_{c,d}(x,y)=x+y\qquad(x,y\in\R_+)
}
is satisfied if and only if one of the following possibilities hold:
\begin{enumerate}[(i)]
\item $\{a,b\}=\{c,d\}=\{1,0\}$, i.e., the two means are equal to the arithmetic mean,
\item there exists $v\in\R$ such that $\{a,b\}=\{1+v,v\}$, $\{c,d\}=\{1-v,-v\}$,
\item $\{a,b\}=\{\frac32,\frac12\}$ and $c+d=0$ (in this case, $G_{c,d}$ is the geometric mean),
\item $a+b=0$ and $\{c,d\}=\{\frac32,\frac12\}$ (in this case, $G_{a,b}$ is the geometric mean).
\end{enumerate}
\end{corollaryG}

The approach followed in \cite{BajPal09b}
heavily used the Maple computer-algebra system to perform the tedious computations of the
various partial derivatives of Gini means up to 12th order.

The main result of this paper completely solves the invariance equation in the class of 
Stolarsky means. 

\begin{theoremS}
Let $a,b,c,d,p,q\in\R$. Then the invariance equation \eq{IE} is valid if and only if one
of the following possibilities holds:
\begin{enumerate}[(i)]
	\item $a+b=c+d=p+q=0$, i.e., all the three means are equal to the geometric mean,
	\item $\{a,b\}=\{c,d\}=\{p,q\}$, i.e., all the three means are equal to each other,
  \item $\{a,b\}=\{-c,-d\}$ and $p+q=0$, i.e., $S_{p,q}$ is the geometric mean and 
              $S_{a,b}=S_{-c,-d}$.
\end{enumerate}
\end{theoremS} 

As a consequence, we obtain the following solution for the
Matkowski--Sut\^o equation, i.e., when $S_{p,q}$ is equal to the arithmetic 
mean in the invariance equation \eq{IE},
which happens exactly when $\{p,q\}=\{2,1\}$.

\begin{corollaryS}
Let $a,b,c,d\in\R$. Then the Matkowski--Sut\^o equation 
\Eq{*}{
  S_{a,b}(x,y)+S_{c,d}(x,y)=x+y\qquad(x,y\in\R_+)
}
holds if and only if $\{a,b\}=\{c,d\}=\{2,1\}$, i.e., both means are equal 
to the arithmetic mean.
\end{corollaryS}

It is interesting to observe here that the parameter sets when \eq{IEG} holds is
much bigger then the corresponding set for \eq{IE}.

In the proof of the above results first we construct an approximation of the Stolarsky means.
Then we use the computer algebra system Maple V Release 9 to 
compute the Taylor expansion of a function (which is in terms of the approximated means)
up to 12th order. Finally, this enables us to describe all the cases of the equality. 

\section{The proof of Theorem S}

First we recall the characterization of the equality of two variable Stolarsky means.

\Lem{Stol}{(Cf.\ \cite{Pal88b}) Let $a,b,c,d\in\R$. Then the identity
\Eq{G}{
  S_{a,b}(x,y)=S_{c,d}(x,y)\qquad(x,y\in\R_+)
}
holds if and only if one of the following possibilities is valid:
\begin{enumerate}[(i)]
\item $a+b=c+d=0$ and, in this case, the two means are equal to the 
geometric mean,
\item $\{a,b\}=\{c,d\}$.
\end{enumerate}
}

The calculations of the higher-order partial derivatives of Stolarsky means at the point 
$(x,y)=(1,1)$ is too complicated (even for computer-algebra systems) because the main 
expression that defines these means is singular at diagonal points.
Therefore, to simplify these computations, we will approximate the Stolarsky mean up to
a sufficiently high order. First, we express the Stolarsky mean in another form.

\Lem{L}{
	For $p,q\in\R$ and $x,y\in\R_+$, the Stolarsky mean $S_{p,q}(x,y)$ can be rewritten as
	\Eq{S}{
		S_{p,q}(x,y)
  :=\begin{cases}
   \left(\dfrac{L(p\log x,p\log y)}{L(q\log x,q\log y)}\right)^\frac1{p-q}
   \quad&\mbox{if } p\neq q,\\[4mm]
    \exp\left(\dfrac{\partial_1 L(p\log x,p\log y)\log x
            +\partial_2 L(p\log x,p\log y)\log y}{L(p\log x,p\log y)}\right)
    \quad&\mbox{if } p=q,
    \end{cases}
	}
	where the function $L:\R^2\to\R$ is defined by
	\Eq{L}{
		L(u,v):=\sum\limits_{n=1}^\infty \frac{u^{n-1}+u^{n-2}v+\cdots+uv^{n-2}+v^{n-1}}{n\,!}\,.
	}
}

\begin{proof}
	Utilizing the Taylor series expansion and substituting $u:=\log x$ and $v:=\log y$, 
we have, for $p\neq 0$ and $x\neq y$,
\Eq{*}{
	\frac{x^p-y^p}{p\left(\log x-\log y\right)}
           &=\frac{\operatorname{e}^{pu}-\operatorname{e}^{pv}}{p\left(u-v\right)}
	   =\frac1{p\left(u-v\right)}\left(1+\frac{pu}{1\,!}+\frac{(pu)^2}{2\,!}+\cdots-
		  1-\frac{pv}{1\,!}-\frac{(pv)^2}{2\,!}-\cdots\right)\\
	   &=\frac1{u-v}\left(\frac{u -v}{1\,!}+\frac{p(u^2-v^2)}{2\,!} 
                 +\frac{p^2(u^3-v^3)}{3\,!}+\cdots\right)\\ 
           &=1+\frac{p(u+v)}{2\,!}+\frac{p^2(u^2+uv+v^2)}{3\,!}+\cdots
            =L(pu,pv).
}
Similarly, we can get that $L(pu,pu)=\operatorname{e}^{pu}$, hence
	\[
	L(p\log x,p\log y):=\begin{cases}
	 \dfrac{x^p-y^p}{p\left(\log x-\log y\right)}, & \mbox{if } p\neq 0, x\neq y,\\[3mm]
	x^p, & \mbox{if } x=y,\\
	1, & \mbox{if } p=0.
	\end{cases}
	\]
Using the above identity and the definition of the Stolarsky means, equality \eq{S} follows
immediately for $p\neq q$. The formula for the case $p=q$ is derived by taking the limit
$q\to p$.
\end{proof}

In order to obtain high-order approximation of Stolarsky means, for $k\in\N$, define
\Eq{Lk}{
	L_k(u,v):=\sum\limits_{n=1}^k \frac{u^{n-1}+u^{n-2}v+\cdots+uv^{n-2}+v^{n-1}}{n\,!}\,
     \qquad(u,v\in\R),
}
and, for $x,y\in\R_+,\,p,q\in\R$,
\Eq{Sk}{
  S_{p,q}^k(x,y)
  :=\begin{cases}
   \left(\dfrac{L_k(p\log x,p\log y)}{L_k(q\log x,q\log y)}\right)^\frac1{p-q}
   \quad&\mbox{if } p\neq q,\\[4mm]
    \exp\left(\dfrac{\partial_1 L_k(p\log x,p\log y)\log x
            +\partial_2 L_k(p\log x,p\log y)\log y}{L_k(p\log x,p\log y)}\right)
    \quad&\mbox{if } p=q.
    \end{cases}
}

\Lem{Der}{For $p,q\in\R$, $k\in\N$,
\Eq{ij}{
  \partial_1^i\partial_2^jS_{p,q}(1,1)=\partial_1^i\partial_2^jS_{p,q}^k(1,1)
}
if $i\geq1$, $j\geq 1$ and $i+j\leq k$.
}

\begin{proof} Observe that, for $\alpha\geq0$, $\beta\geq 0$ with 
$\alpha+\beta\leq k$, we have
\Eq{*}{
  \partial_1^\alpha\partial_2^\beta L(0,0)=\partial_1^\alpha\partial_2^\beta L_k(0,0).
}
In the case $p\neq q$,
by applying the elementary rules of differentiation and the identities \eq{S} and \eq{Sk}, 
we get that the formulae for the partial derivatives $\partial_1^i\partial_2^jS_{p,q}(1,1)$ and 
$\partial_1^i\partial_2^jS_{p,q}^k(1,1)$ are the same expressions of all the partial derivatives 
$\partial_1^\alpha\partial_2^\beta L(0,0)$ and $\partial_1^\alpha\partial_2^\beta L_k(0,0)$,
respectively, where $0\leq\alpha\leq i$ and $0\leq\beta\leq j$. Hence, equality \eq{ij} follows 
for $p\neq q$. In the case $p=q$, the identity is derived by taking the limit $q\to p$.
\end{proof}

To make the calculations even more simple, we define the function $E:\R^4\to \R$, by
\Eq{*}{
	E(p,q,u,v):=\log\big(S_{p,q}(e^u,e^v)\big)=
	\begin{cases}
          \dfrac{\log L(pu,pv)-\log L(qu,qv)}{p-q} \ &\mbox{if } p\neq q, \\[5mm]
          \dfrac{\partial_1 L(pu,pv)u+\partial_2 L(pu,pv)v}{L(pu,pv)} \ &\mbox{if } p=q.
	\end{cases}
}
The invariance equation \eq{IE} can be rewritten in the form
\[
	E\big(p,q,E(a,b,u,v),E(c,d,u,v)\big)=E(p,q,u,v)\qquad(u,v\in\R),
\]
hence
\[
	F(x):=E\big(p,q,E(a,b,x,-x),E(c,d,x,-x)\big)-E(p,q,x,-x)=0\qquad(x\in\R).
\]
Thus, the derivatives $F^{(m)}$ vanish at $x=0$ for all $m\geq 0$, $m\in \mathbb Z$.
By the symmetry of Stolarsky means, $F$ is an even function.
Therefore, $F^{(2m+1)}(0)=0$ holds automatically. In order to obtain the necessity
of the conditions of the theorem, we will need to investigate the equalities
\Eq{Fder}{
  F^{(2m)}(0)=0, \qquad (1\leq m\leq 6).
}
To compute these derivatives, define the function $F_k:\R\to \R,\, k\in \N$, as
\[
	F_k(x):=E_k\big(p,q,E_k(a,b,x,-x),E_k(c,d,x,-x)\big)-E_k(p,q,x,-x),
\]
where 
\[
	E_k(p,q,u,v):=\log\big(S^k_{p,q}(e^u,e^v)\big)=
	\begin{cases}
          \dfrac{\log L_k(pu,pv)-\log L_k(qu,qv)}{p-q} \ &\mbox{if } p\neq q, \\[5mm]
          \dfrac{\partial_1 L_k(pu,pv)u+\partial_2 L_k(pu,pv)v}{L_k(pu,pv)} \ &\mbox{if } p=q.
	\end{cases}
\]
It follows from \lem{Der} that the derivatives of the functions $E_k$ and $E$
at $(p,q,0,0)$, and hence the derivatives of $F_k$ and $F$ up to the order $k$ 
at the point $x=0$ coincide.
Therefore, to analyze the equations \eq{Fder}, it is sufficient to consider
the identities
\Eq{Fkder}{
	F_k^{(2m)}(0)=0, \qquad (1\leq m\leq 6),
}
where $2m \leq k$.
Thus in the proof of the theorem, we compute the Taylor expansion of the appropriate
$F_k$ at $x=0$, hence we can obtain conditions for the
unknown parameters $a,b,c,d,p,q$.

\begin{proof}[The proof of Theorem S]

Assume that \eq{IE} holds.
In the syntax of the Maple language, we define the functions $L_k,E_k$ and $F_k$ in the following way:
\begin{verbatim}
	L:=(u,v,k)->add((1/n!)*add(u^(n-i)*v^(i-1),i=1..n),n=1..k);
	E:=(p,q,u,v,k)->(1/(p-q))*(ln(L(pu,pv,k))-ln(L(qu,qv,k)));
	F:=(x,k)->E(p,q,E(a,b,x,-x,k),E(c,d,x,-x,k),k)-E(p,q,x,-x,k);
\end{verbatim}
This produces the following output:
\Eq{*}{
	L:=&(u,v,k) \to \operatorname{add}\left( \frac{\operatorname{add}
			(u^{n-i}v^{i-1},i=1\,.\,.\,n)}{n\, !}, n=1\,.\,.\,k\right),\\
	E:=&(p,q,u,v,k) \to \frac{\log \big(L(pu,pv,k)\big)-\log \big(L(qu,qv,k)\big)}{p-q},\\
	F:=&(x,k) \to E(p,q,E(a,b,x,-x,k),E(c,d,x,-x,k),k)-E(p,q,x,-x,k).
}

First we evaluate the second-order Taylor coefficient $C_2$ of $F$ at $x=0$ by replacing
$F$ by $F_3$:
\begin{verbatim}
	> C[2]:=simplify(coeftayl(F(x,3),x=0,2));
\end{verbatim}
\[
	C_2=\frac1{12}\,d-\frac16\, q + \frac1{12}\,p-\frac16\, q+ \frac1{12}\,a + \frac1{12}\,c.
\]
(Note that the Maple-definition of the function $E$ is valid only if
$(p-q)(a-b)(c-d)\neq0$, however, the Taylor coefficient $C_2$ and also the subsequent ones,
are correct also in the singular case $(p-q)(a-b)(c-d)=0$.)

By the invariance equation $C_2=0$, therefore
\[
	\frac{a+b+c+d}4=\frac{p+q}2.
\]
In order to simplify the evaluation of the higher-order Taylor coefficients,
we introduce the notations
\Eq{vt}{
  w:=&\frac{a+b+c+d}{4}=\frac{p+q}2, \\
  v:=&\frac{a+b-(c+d)}{4},\\
  t:=&\Big( \frac{p-q}2\Big)^2,\\
  r:=&\frac{(a-b)^2+(c-d)^2}{8},\\
  s:=&\frac{(a-b)^2-(c-d)^2}{8}.
}
(In the definition of $w$ we utilized the condition $C_2=0$.)
Then we can express the parameters $a,b,c,d,p,q$ in the following form:
\begin{verbatim}
> a:=w+v+sqrt(r+s); b:=w+v-sqrt(r+s); 
  c:=w-v+sqrt(r-s); d:=w-v-sqrt(r-s); 
  p:=w+sqrt(t); q:=w-sqrt(t);
\end{verbatim}
\Eq{*}{
  a&:= w+v+\sqrt {r+s}\\  b&:= w+v-\sqrt {r+s}\\
  c&:= w-v+\sqrt {r-s}\\  d&:= w-v-\sqrt {r-s}\\
}
\Eq{*}{
  p&:= w+\sqrt {t}\\  q&:= w-\sqrt {t}
}
Now we evaluate the 4th order Taylor coefficient $C_4$ of $F$ (i.e., of $F_5$)
at $x=0$ by inputting: 
\begin{verbatim}
> C[4]:=simplify(coeftayl(F(x,5),x=0,4));
\end{verbatim}
\[
  C_4:=-\frac{1}{45}vs-\frac{4}{135}wv^2+\frac{1}{45}wt-\frac{1}{45}wr
\]
The condition $C_4=0$ yields that $wt=wr+vs+\dfrac43wv^2$.

If $w=0$, then $p+q=0$, which means that $S_{p,q}$ is the geometric mean.
Therefore, the invariance equation can be written as
\[
	S_{a,b}(x,y)S_{c,d}(x,y)=xy \qquad(x,y\in\R_+).
\]
This results
\Eq{*}{
  S_{a,b}(x,y)=\frac{1}{S_{c,d}(1/x,1/y)}=S_{-c,-d}(x,y)\qquad(x,y\in\R_+).
}
Using \lem{Stol}, this identity yields that either 
$a+b=c+d=0$ or $\{a,b\}=\{-c,-d\}$ must hold.
In this case we get that one of the conditions (i) or (iii) of our theorem is valid.
Conversely, if conditions (i) or (iii) hold then \eq{IE} can easily be seen.

In the rest of the proof, we may assume that $w$ is not zero. Then, 
from condition $C_4=0$, we can express $t$ in terms of $w,v,r,s$:
\begin{verbatim}
> t:=r+v*s/w+(4/3)*v^2;
\end{verbatim}
\Eq{t}{
  t:= r+\frac{vs}w+\frac43\, v^2.
}
Next, we evaluate the 6th order Taylor coefficient $C_6$ of $F$ (i.e., of $F_7$) at $x=0$:
\begin{verbatim}
> C[6]:=simplify(coeftayl(F(x,7),x=0,6));
\end{verbatim}
We get that
\Eq{*}{
  C_6:=\frac{2(9w^2s^2+8w^2v^4+45w^2rv^2+39w^3sv+6wv^3s-13w^4v^2-9v^2s^2)}{8505w}
}

If $v=0$, then $C_6=0$ implies that $s=0$. Hence, from $\eq{t}$ it follows
that $t=r$ and we get that $a=c=p$ and $b=d=q$, i.e., condition (ii) of our
theorem holds. Conversely, if condition (ii) holds then the invariance equation 
\eq{IE} is trivially valid.

In the rest of the proof, we may assume that $v$ is also not zero. 
Observe that the 6th order coefficient $C_6$ does not involve higher-order 
powers of $r$. Therefore, the equation $C_6=0$ can be solved for $r$. Temporarily, 
we denote this solution by $R$:
\begin{verbatim}
> R:=(13*w^4*v^2-9*w^2*s^2-8*w^2v^4+9*v^2*s^2-39*w^3*v*s-6*w*v^3*s)/
		 (45*w^2*v^2);
\end{verbatim}
\Eq{R}{
  R:= \frac{13w^4v^2-9w^2s^2-8w^2v^4+9v^2s^2-39w^3vs-6wv^3s}{45w^2v^2}
}

Finally, we evaluate the 13th order Taylor polynomial of $F_{13}$ at $x=0$ (the 
Maple output is suppressed by putting : instead of ; to the end of the Maple command, 
for the sake of brevity), then we extract the 8th, 10th and 12th order Taylor 
coefficients, denoted by $C_8$, $C_{10}$ and $C_{12}$, respectively,
and replace $r$ by $R$ by inputting:
\begin{verbatim}
> T:=simplify(taylor(F(x,13),x=0,13)):
  for i from 8 to 12 by 2 
      do C[i]:=simplify(subs(r=R,simplify(coeff(T,x,i))),factor) od;
\end{verbatim}
\Eq{*}{
C_8&:=\frac {1}{9568125w^3v^2} 
\big(4347w^5vs^3-3616w^8v^4+1242w^2v^6s^2+16137w^6v^2s^2+224w^4v^8\\
&\hspace{1cm}-4644w^3v^3s^3+756w^3v^7s+81w^4s^4+11976w^7v^3s+4992w^6v^6-162w^2v^2s^4\\
&\hspace{1cm}-17379w^4v^4s^2-12732w^5v^5s+297wv^5s^3+81v^4s^4\big)
}
\Eq{*}{
C_{10}&:=\frac {2}{2841733125w^5v^4}
\big(141632w^{10}v^8-2187w^2v^4s^6-58806w^3v^5s^5+107406w^5v^3s^5\\
&\hspace{1cm}-338432w^{12}v^6++34928w^6v^{12}-523908w^9v^3s^3+2187w^4v^2s^6+126459w^4v^{10}s^2\\
&\hspace{1cm}+316272w^8v^{10}+143541w^8v^6s^2+92508w^5v^{11}s+24948w^2v^8s^4+1319016w^{11}v^5s\\
&\hspace{1cm}-847044w^9v^7s+927324w^{10}v^4s^2-1197324w^6v^8s^2-564480w^7v^9s-729w^6s^6\\
&\hspace{1cm}+55026w^3v^9s^3-794610w^5v^7s^3+672786w^6v^4s^4+3402wv^7s^5-52002w^7vs^5\\
&\hspace{1cm}-313065w^4v^6s^4+1263492w^7v^5s^3-384669w^8v^2s^4+729v^6s^6\big)
}
\Eq{*}{
C_{12}&:=\frac {2}{872767286015625w^7v^6} 
\big(54403812w^4v^4s^8-36269208w^6v^2s^8-2472351012w^7v^3s^7\\
&\hspace{4mm}-36269208w^2v^6s^8-6728400999w^4v^8s^6-23518313469w^8v^4s^6+28808296272w^{10}v^{14}\\
&\hspace{4mm}-43502512128w^{16}v^8-32563674368w^{14}v^{10}+55532868672w^{12}v^{12}+4082581552w^8v^{16}\\
&\hspace{4mm}+51381378wv^9s^7+1616450580w^3v^{11}s^5+66139190820w^7v^7s^5+5633126127w^4v^{12}s^4\\
&\hspace{4mm}-27686800350w^5v^9s^5+23171026860w^{11}v^3s^5-198268186320w^{13}v^5s^3+9067302w^8s^8\\
&\hspace{4mm}-50400222003w^{12}v^4s^4+8922828051w^{10}v^2s^6+527871816w^2v^{10}s^6+2575113768w^5v^5s^7\\
&\hspace{4mm}+33853773582w^{10}v^6s^4+8610065640w^5v^{13}s^3+78734957472w^8v^8s^4+806989878w^9vs^7\\
&\hspace{4mm}-56327456871w^8v^{12}s^2-297179050266w^{10}v^{10}s^2-14684723544w^{13}v^9s\\
&\hspace{4mm}+203758408416w^{15}v^7s-195782704164w^{11}v^{11}s-67821635178w^6v^{10}s^4\\
&\hspace{4mm}+10150547496w^7v^{15}s+13660568109w^6v^{14}s^2-84844566450w^7v^{11}s^3\\
&\hspace{4mm}+337457994084w^{12}v^8s^2+2387944944w^{14}v^6s^2-3441528204w^9v^{13}s\\
&\hspace{4mm}-67735201590w^9v^9s^3+342237888720w^{11}v^7s^3-961134012w^3v^7s^7\\
&\hspace{4mm}+20796014601w^6v^6s^6-63239867910w^9v^5s^5+9067302v^8s^8\big)
}
The Taylor coefficients $C_8$, $C_{10}$, and $C_{12}$ are of the form
\Eq{*}{
  C_8=\frac{1}{9568125w^3v^2}P_8,\quad 
  C_{10}=\frac{2}{2841733125w^5v^4}P_{10},\quad
  C_{12}=\frac {2}{872767286015625w^7v^6}P_{12},
}
where $P_8$, $P_{10}$, and $P_{12}$ are polynomials of the variables $v,w,s$.
They can be obtained by the following Maple commands (whose output is suppressed):
\begin{verbatim}
>  P[8]:=op(2,C[8]): P[10]:=op(2,C[10]): P[12]:=op(2,C[12]):
\end{verbatim}
The equalities $C_8=C_{10}=C_{12}=0$ imply that 
$P_8=P_{10}=P_{12}=0$.

The variable $s$ is a common root of the polynomials $P_8$ and $P_{10}$. 
Therefore the resultant $R_{8,10}$ of these two polynomials (with respect to $s$) is zero:
\begin{verbatim}
> R[8,10]:=factor(resultant(P[8],P[10],s));
\end{verbatim}
{\small\Eq{*}{
&R_{{8,10}}:=28242953648100000000 w^{24} v^{24} (v-w)^8 (v+w)^8 \big(395726752304\,{v}^{32}-28019198519832\,{w}^{2}{v}^{30}\\
	&+1192972799035666\,{w}^{4}{v}^{28}-36617671790074251\,{w}^{6}{v}^{26}+601554420387156651\,{w}^{8}{v}^{24}\\
	&-3652037976710860175\,{w}^{10}{v}^{22}-1101310194408221307\,{w}^{12}{v}^{20}+62048533824813847173\,{w}^{14}{v}^{18}\\
  &-175575191501013599783\,{w}^{16}{v}^{16}+52614376847529172973\,{w}^{18}{v}^{14}+435211540238087039223\,{w}^{20}{v}^{12}\\
  &-793895884964266327270\,{w}^{22}{v}^{10}-773252618095825970136\,{w}^{24}{v}^{8}-492199682627262911866\,{w}^{26}{v}^{6}\\
  &+183522699320559043726\,{w}^{28}{v}^{4}-43030934088053846752\,{w}^{30}{v}^{2}+912066926976343384\,{w}^{32}\big)
}}
The resultant is zero if either $vw(v-w)(v+w)=0$ holds or $v$ and $w$ are solutions of a 
homogeneous two variable polynomial equation of degree 32. Writing $w$ in the form
\begin{verbatim}
> w:=z*v;
\end{verbatim}
\Eq{*}{
  w:= zv
}
we get that $z$ is a root of a 32nd degree polynomial $P_{8,10}$, where:
\begin{verbatim}
> P[8,10]:=simplify(op(4,R[8,10])/v^32);
\end{verbatim}
\Eq{*}{
P&_{8,10}:=395726752304-28019198519832{z}^{2}+1192972799035666{z}^{4}-36617671790074251{z}^{6}\\
	&+601554420387156651{z}^{8}-3652037976710860175{z}^{10}-1101310194408221307{z}^{12}\\
	&+62048533824813847173{z}^{14}-175575191501013599783{z}^{16}+52614376847529172973{z}^{18}\\
	&-435211540238087039223{z}^{20}-793895884964266327270{z}^{22}+773252618095825970136{z}^{24}\\
	&-492199682627262911866{z}^{26}+183522699320559043726{z}^{28}-43030934088053846752{z}^{30}\\
	&+912066926976343384{z}^{32}.
}
The variable $s$ is also a common root of the two polynomials $P_8$ and $P_{12}$. 
Therefore the resultant $R_{8,12}$ of these polynomials (with respect to $s$)
is again zero:
\begin{verbatim}
> R[8,12]:=factor(resultant(P[8],P[12],s)):
\end{verbatim}
We now get that $v$ and $w$ are solutions of a 
homogeneous two variable polynomial of degree 44, whence we get that $z$ is a 
root of the 44th degree polynomial $P_{8,12}$, where:
\begin{verbatim}
> P[8,12]:=op(4,R[8,12]);
\end{verbatim}
Now computing the resultant of the two polynomials $P_{8,10}$ and $P_{8,12}$ by
\begin{verbatim}
> Q:=resultant(P[8,10],P[8,12],z):
\end{verbatim}
it follows that $Q$ is a (huge) nonzero number, hence $P_{8,10}$ and $P_{8,12}$ cannot 
have a common root. This proves that $R_{8,10}$ and $R_{8,12}$ can be simultaneously zero
if and only if $vw(v-w)(v+w)=0$ holds. Hence $C_8=C_{10}=C_{12}=0$ can hold only
in this case.

Since we have $vw\neq0$, hence $(v-w)(v+w)=0$ must hold, i.e, $v=\pm w$.
Thus, from \eq{R} and \eq{t}, we get that $r=\dfrac{w^2}9\mp s$ and $t=\dfrac{13}9 w^2$,
respectively.

In the case when $v=w$, the equations in \eq{vt} yield that
\Eq{C1}{
	c=-d, \qquad a=\dfrac{7}{3}\, w, \qquad  b=\dfrac{5}{3}\, w, \qquad
	p=\bigg(1+\dfrac{\sqrt{13}}{3}\bigg) w,  \qquad q=\bigg(1-\dfrac{\sqrt{13}}{3}\bigg) w.
}
The first equality yields that $S_{c,d}$ is the geometric mean. Hence, we may assume that
$c=-d=1$. To simplify the computations, we can also assume that $w=3$.
We show that these parameters are not solutions of the invariance equation.
For $k\in \N$, we now have that
\[
	F_k(x):=E_k\big(3+\sqrt{13},3-\sqrt{13},E_k(7,5,x,-x),E_k(1,-1,x,-x)\big)
	-E_k\big(3+\sqrt{13},3-\sqrt{13},x,-x\big).
\]
In Maple, we input
\begin{verbatim}
> F:=(x,k)->E(3+sqrt(13),3-sqrt(13),E(7,5,x,-x,k),E(1,-1,x,-x,k),k)-
  E(3+sqrt(13),3-sqrt(13),x,-x,k);
\end{verbatim}
We compute the 10th order Taylor coefficient of $F_{11}$ by
\begin{verbatim}
> simplify(coeftayl(F(x,11),x=0,10));
\end{verbatim}
whence we get that this coefficient is $-\dfrac{12352}{5775}$, i.e., it is not
zero, which means that the parameters in \eq{C1} do not provide solution to the
invariance equation.

In the case when $v=-w$, we have that
\Eq{*}{
	a=-b, \qquad c=\dfrac{7}{3}\, w, \qquad  d=\dfrac{5}{3}\, w, \qquad
	p=\left(1+\dfrac{\sqrt{13}}{3}\right) w,  \qquad q=\left(1-\dfrac{\sqrt{13}}{3}\right) w,
}
whence a similar calculation as in the previous case shows that we again do not 
get an additional solution to the invariance equation.
\end{proof}

\textbf{Acknowledgement.} The authors thank the anonymous referees for their help
and suggestions.

\end{document}